\documentclass[11pt]{amsart}
\usepackage{etex}
\usepackage{diagbox}
\usepackage{array}
\usepackage{graphicx, pictex}
\usepackage{stmaryrd}
\usepackage{float}
\restylefloat{table}
\usepackage{multirow}
\usepackage{amsmath,amssymb,amsthm,mathrsfs}
\allowdisplaybreaks
\def\d{\Omega}

\def\du#1#2#3{\overset{#3}{\underset{#2}{#1}}}
\def\Forall{\quad \hbox{ for all }}

\def\M{{\mathcal{M}}}

\newcommand{\tn}[1]{\lVert\kern-1pt\lvert{#1}\rvert\kern-1pt\rVert}
\def\<{{\langle}}
\def\>{{\rangle}}
\def\Forall{\quad \hbox{ for all }}

\def\A{{\mathcal{A}}}
\def\C{{\mathcal{C}}}

\def\d{\Omega}

\def\Forall{\quad \hbox{ for all }}

\def\d{\Omega}

\def\Forall{\quad \hbox{ for all }}

\def\tb#1{{\|\kern-1pt| #1 \|\kern-1pt|}}
\def\nm2#1#2{\|#1\|_{2,\d_{#2}}}

\def\R{\mathbb{R}}

\def\L{\mathcal{L}}

 \theoremstyle{plain}
 \newtheorem{thm}{Theorem}[section]
 \numberwithin{equation}{section} %% Comment out for sequentially-numbered
 \numberwithin{figure}{section} %% Comment out for sequentially-numbered
 \theoremstyle{plain}
 \newtheorem{prop}[thm]{Proposition}
 \theoremstyle{plain}
 \newtheorem{algorithm}[thm]{Algorithm} %%Delete [thm] to re-start numbering
 \theoremstyle{plain}
 \newtheorem{theorem}[thm]{Theorem}
 \theoremstyle{plain}
 
\theoremstyle{plain}
 \newtheorem{remark}[thm]{Remark}
 \theoremstyle{plain}
 
\def\A{{\mathcal{A}}}

\def\C{{\mathcal{C}}}

\def\M{{\mathcal{M}}}

\def\L{{\mathcal{L}}}
\def\d{{\Omega}}

\def\Forall{\quad \hbox{ for all }}
\def\<{{\langle}}
\def\>{{\rangle}}
\def\R{\mathbb{R}}

\def\A{{\mathcal{A}}}
\def\du#1#2#3{\overset{#3}{\underset{#2}{#1}}}

\usepackage{amssymb}

\begin{document}
\title[Least Squares Preconditioning of Mixed Formulations]
{Saddle Point Least Squares  Preconditioning of Mixed Methods}

\author{Constantin Bacuta}
\address{University of Delaware,
Department of Mathematics,
501 Ewing Hall 19716}
\email{bacuta@udel.edu}

\author{Jacob Jacavage}
\address{University of Delaware,
Department of Mathematics,
501 Ewing Hall 19716}
\email{jjacav@udel.edu}


\thanks{The work was supported  by NSF, DMS-1522454.}

\begin{abstract}
We present a simple way to discretize and precondition  mixed variational formulations. Our theory connects with, and takes advantage of, the classical theory of symmetric saddle point problems and the theory of preconditioning symmetric positive definite operators. Efficient iterative processes  for solving the discrete mixed formulations  are proposed and choices for discrete spaces that are always compatible are provided.  For the proposed discrete spaces and solvers, a basis is needed only for the test spaces and assembly of a global saddle point system is avoided. We prove sharp approximation properties for the discretization and iteration errors and also provide a sharp estimate for the convergence rate of the proposed algorithm in terms of the condition number of the elliptic preconditioner and the discrete $\inf-\sup$ and $\sup-\sup$  constants of the pair of discrete spaces.  
\end{abstract}
\maketitle

%%%%%%%%%%%%%%%%%%%%%%%%%%%%%%%%%%%%%%%%%%5
\section{Introduction}\label{introduction}

We  provide  a general  approach  in preconditioning  mixed problems of the form: Given $f\in V^*$, find $p\in Q$ such that
\begin{equation}\label{cont_problem}
b(v,p)=\langle f,v\rangle \Forall v\in V,
\end{equation}
where $V$ and $Q$ are Hilbert spaces and $b(\cdot,\cdot)$ is a continuous bilinear form on $V\times Q$  satisfying
an $\inf-\sup$ condition. In \cite{BM12, Dahmen-Welper-Cohen12}, a connection was made between problems of the form \eqref{cont_problem} and  a natural  saddle point formulation. More specifically, if  $a(\cdot, \cdot)$ is the inner product on $V$, then $p$ is the unique solution of \eqref{cont_problem} if and only if $(u=0,p)$ is the unique solution to: Find $(u,p)\in V\times Q$ such that
\begin{equation}\label{abstract:variationalSPP}
\begin{array}{lclll}
a(u,v) & + & b( v, p) &= \langle f,v \rangle &\ \Forall  v \in V,\\
b(u,q) & & & =0   &\  \Forall  q \in Q,  
\end{array}
\end{equation}
where appropriate assumptions on $f$ and the form $b(\cdot,\cdot)$ hold, see Section \ref{notation}. Thus,  \eqref{abstract:variationalSPP} is a  saddle point reformulation of \eqref{cont_problem}. It is clear that the solution component $p$ of  this reformulation is independent of the inner product (hence the norm) considered on $V$. This observation is essential for the discretization and preconditioning of \eqref{cont_problem}.  

For finite dimensional approximation spaces $V_h\subset V$ and $\M_h\subset Q$ satisfying
a discrete  $\inf-\sup$ condition, we consider the discrete problem of finding $(u_h,p_h) \in V_h\times \M_h$ such that
\begin{equation}\label{discrete:variationalSPP}
\begin{array}{lclll}
a(u_h,v_h) & + & b( v_h, p_h) &= \langle f,v_h \rangle &\ \Forall  v_h \in V_h,\\
b(u_h,q_h) & & & =0   &\  \Forall  q_h \in \M_h, 
\end{array}
\end{equation} 
which approximates the solution  $(u=0,p)$ of \eqref{abstract:variationalSPP}.  The discrete  variational formulation \eqref{discrete:variationalSPP}  is  in fact  a {\it saddle point least squares discretization}  of \eqref{cont_problem}, see Section \ref{the_discrete_problem}. This saddle point discretization of \eqref{cont_problem} is also adapted by Demkowicz and Gopalakrishnan in \cite{demkowicz-gopalakrishnanDPG10,demkowicz-gopalakrishnanDPG13}. When solving the above problem, finding bases for the discrete trial space $\M_h$ and assembling a block stiffness matrix for  \eqref{discrete:variationalSPP} can be avoided by applying an Uzawa type algorithm. However, any attempt to solve \eqref{discrete:variationalSPP} by an Uzawa iterative process requires the exact  inversion of the  operator  $A_h$ associated with the inner product $a(\cdot, \cdot)$ on $V_h$. To avoid exact inversion and to speed up the iterative solvers, we consider another form  $\tilde{a}(\cdot,\cdot)$ on $V_h$, which leads to an equivalent norm on $V_h$, and  introduce a preconditioned discrete saddle point problem: Find $(\tilde{u}_h,\tilde{p}_h)\in V_h\times \M_h$ such that 
\begin{equation}\label{discrete:variationalSPPprec}
\begin{array}{lclll}
\tilde{a}(\tilde{u}_h,v_h) & + & b( v_h, \tilde{p}_h) &= \langle f,v_h \rangle &\ \Forall  v_h \in V_h,\\
b(\tilde{u}_h,q_h) & & & =0   &\  \Forall  q_h \in \M_h, 
\end{array}
\end{equation}
where the  action of the operator  $\tilde {A}_h^{-1}$ associated with the inner product $\tilde{a}(\cdot, \cdot)$ on  $V_h$  is assumed to be fast and  easy to implement.

The goals of this paper are: describe how well the component solution $\tilde{p}_h$ of \eqref{discrete:variationalSPPprec} approximates the solution  $p$ of \eqref{cont_problem}, describe possible choices for the discrete pairs $(V_h, \M_h)$, and propose an efficient iterative solver for \eqref{discrete:variationalSPPprec} and estimate its convergence rate.  

The paper is organized as follows. In Section 2, an abstract theory for saddle point least squares formulations is presented. Section 3 describes the general preconditioning theory and approximation results. In addition, convergence rates for the proposed iterative solver are estimated. Possible choices for discrete pairs of spaces are discussed in Section 4.

%%%%%%%%%%%%%%%%%%%%%%%%%%%%%%%%%%%%%%%%%%%%%%%%
\section{Abstract Saddle Point Least Squares Formulation for Mixed Methods}\label{abstract}
\subsection{Notation and the continuous problem} \label{notation}

Let $V$ and $Q$ be infinite dimensional Hilbert spaces and assume the inner products $a(\cdot, \cdot)$ and $(\cdot, \cdot)$ induce the norms $|\cdot|_V =|\cdot| =a(\cdot, \cdot)^{1/2}$ and $\|\cdot\|_{Q}=\|\cdot\|=(\cdot, \cdot)^{1/2}$. The duals of $V$ and $Q$ will be denoted by $V^*$ and $Q^*$, respectively. The dual pairings on $V^* \times V$ and $Q^* \times Q$ will both be denoted by $\langle \cdot, \cdot \rangle$. With the inner products $a(\cdot, \cdot)$ and $(\cdot, \cdot)$, we associate the operators $\A: V\to V^*$ and $\C: Q\to Q^*$ defined by 
\[
\langle \A u,v \rangle=a(u, v)  \Forall u ,v \in V,
\]
and
\[
\langle \C p,q\rangle=(p, q) \Forall  \ p, q \in Q.
\]

Assume that $b(\cdot, \cdot)$ is a continuous bilinear form on $V\times Q$ satisfying
 the $\inf-\sup$ condition
 \begin{equation} \label{inf-sup_a}
\du{\inf}{p \in Q}{} \ \du {\sup} {v \in V}{} \ \frac {b(v, p)}{|v|\,\|p\|} =m>0,
\end{equation} 
and is bounded, i.e.,
\begin{equation}\label{sup-sup_a}
\du{\sup}{p \in Q}{} \ \du {\sup} {v \in V}{} \ \frac {b(v, p)}{|v|\,\|p\|} =M <\infty.
\end{equation}
With the form $b(\cdot, \cdot)$, we associate the linear operators $B:V\to Q^*$ and $B^*:Q\to V^*$ defined through the duality pairings
\[
\<B v,q\>=b(v, q)= \langle B^*q, v \rangle  \Forall v \in V, \ q \in Q.
\]

It is well known that if, in addition to the assumptions on $b(\cdot,\cdot)$, $f$ satisfies the compatibility condition
\begin{equation}\label{comp_condition}
\langle f,v\rangle=0 \Forall v\in V_0:=\{v\in V\,|\,b(v,q) =0, \ \mathrm{for \ all}\  q \in Q\},
\end{equation}
then \eqref{cont_problem} has a unique solution $p$,  see e.g. \cite{A-B,B09}. Furthermore, $(u=0,p)$ is the unique solution of \eqref{abstract:variationalSPP}. 
\begin{remark}
The saddle point problem \eqref{abstract:variationalSPP}  has a unique solution $(u, p)$ regardless of the compatibility condition \eqref{comp_condition}. The operator form of problem \eqref{cont_problem} is equivalent to  finding $p\in Q$ such that
\[
\A^{-1}B^*p=\A^{-1}f, 
\]
and solving for $p$  from \eqref{abstract:variationalSPP} gives 
\[
(\C^{-1}B)(\A^{-1}B^*)p=(\C^{-1}B)\A^{-1}f.
\]
Since $\C^{-1}B$ is the Hilbert transpose of  $\A^{-1}B^*$, we  have that  the $p$ component of the solution of \eqref{abstract:variationalSPP} is the {\it least squares solution} of \eqref{cont_problem}.
\end{remark}
For the rest of this paper we assume that the compatibility condition \eqref{comp_condition} holds, and consequently, problem \eqref{cont_problem} has a unique solution.

\subsection{Saddle point least squares discretization}\label{the_discrete_problem}
Let $V_h\subset V$ and \\ $\M_h\subset Q$ be finite dimensional approximation spaces and $A_h$  be the discrete version of the operator $\A$, i.e.,  $A_h$ satisfies
\[
\langle A_h u_h,v_h\rangle=a(u_h,v_h) \Forall u_h,v_h\in V_h.
\]
We define the discrete operators $B_h:V_h\to \M_h$ and $B_h^*:\M_h\to V_h^*$ by
\[
(B_hv_h,q_h)=b(v_h,q_h)=\langle B_h^*q_h,v_h\rangle \Forall v_h\in V_h,q_h\in \M_h. 
\]
Note that the operator $B_h$ is defined using the inner product on $\M_h$ and not with the duality on $\M^*_h \times \M_h$. Thus, we can define the discrete Schur complement $S_h:\M_h \to \M_h$ as $S_h=B_h\, A_h^{-1} B_h^*$. We further 
 assume the following discrete $\inf-\sup$ condition holds for the pair of spaces $(V_h,\M_h)$:
 \begin{equation} \label{inf-sup_discrete}
\du{\inf}{p_h \in \M_h}{} \ \du {\sup} {v_h \in V_h}{} \ \frac {b(v_h, p_h)}{|v_h|\,\|p_h\|} =m_h>0.
\end{equation} 
It is well known that the spectrum of $S_h$ satisfies $\sigma(S_h) \subset [m_h^2, M_h^2]$, where 
\begin{equation}\label{sup-sup_h}
M_h:= \du{\sup}{p_h \in \M_h}{} \ \du {\sup} {v_h \in V_h}{} \ \frac {b(v_h, p_h)}{|v_h|\,\|p_h\|} \leq M <\infty,
\end{equation}
and that $m_h^2, M_h^2$ are (the extreme) eigenvalues of $S_h$. 
Define 
\[
V_{h,0}:=\{v_h\in V_h\,|\, b(v_h,q_h)=0\Forall q_h\in \M_h\},
\] 
to be the kernel of the discrete operator $B_h$. We define $f_h \in V_h^*$ to be the restriction of $f$ to $V_h$, i.e.,   $\langle f_h, v_h \rangle:=\langle f, v_h \rangle$ for all $v_h \in V_h$.
\begin{remark}\label{well-posed_prop}
In the case $V_{h,0}\subset V_0$, the compatibility condition \eqref{comp_condition} implies the discrete compatibility condition
\[
\langle f,v_h\rangle =0 \Forall v_h\in V_{h,0}.
\]
Hence, under assumption \eqref{inf-sup_discrete}, the problem of finding $p_h\in \M_h$ such that 
\begin{equation}\label{discrete_var_form}
b(v_h,p_h)=\langle f,v_h\rangle,  \ v_h\in V_h, \ \text{or} \ B_h^*\, p_h = f_h, \ \text{or} \ A_h^{-1}B_h^*\, p_h = A_h^{-1}f_h,
\end{equation}
has a unique solution. In general, \eqref{comp_condition} may not hold on $V_{h,0}$ and problem  \eqref{discrete_var_form} may not be well-posed. However, if the form $b(\cdot, \cdot)$ satisfies \eqref{inf-sup_discrete}, then the problem of finding $(u_h,p_h) \in V_h\times \M_h$ satisfying  \eqref{discrete:variationalSPP} does have a unique solution. Solving for $p_h$ from   \eqref{discrete:variationalSPP}, we obtain 
 \begin{equation} \label{Sh}
S_h \, p_h= B_h(A_h^{-1}B_h^*) \,p_h = B_h A_h^{-1}f_h.
\end{equation}
Since the Hilbert transpose of $B_h$ is  $B_h^T = A_h^{-1} B_h^*$,  we call  the component $p_h$  of the solution $(u_h,p_h)$ 
of \eqref{discrete:variationalSPP} the {\it saddle point least squares} approximation of the solution $p$ of the original mixed problem \eqref{cont_problem}.
\end{remark}
The following error estimate for $\|p-p_h\|$ was proved in \cite{BQ15}. 

%Using  the extra compatibility condition  \eqref{eq:BBsuf},   a sharp error estimate for $\|p-p_h\|$ was proved in \cite{BQ15}  based on the  Xu-Zikatanov argument, see \cite{xu-zikatanov-BBtheory}. 
 
\begin{theorem}\label{th:sharpEE} 
Let $b:V \times Q \to \R$  satisfy \eqref{inf-sup_a} and \eqref{sup-sup_a} and assume that   ${f}  \in V^*$  is given and satisfies \eqref{comp_condition}. Assume that  $p$  is the  solution  of \eqref{cont_problem} and  $V_h \subset V$,  $ {\M}_h \subset  Q$ are  chosen such that the discrete $\inf-\sup$ condition   \eqref{inf-sup_discrete} holds. If  $\left (u_h, p_h \right )$ is the  solution  of \eqref{discrete:variationalSPP}, then the following error estimate holds:
\begin{equation}\label{eq:er4LS}
\frac 1 M |u_h| \leq \|p-p_h\| \leq  \frac{M}{m_h} \  \du{\inf}{q_h \in\M_h
}{}  \|p-q_h\|.
\end{equation} 
\end{theorem}

\subsection{An Uzawa CG iterative solver}\label{iter_discrete_problem}
Note that a global linear system may be difficult to assemble when solving \eqref{discrete:variationalSPP} as bases for the trial spaces $\M_h$, which are chosen to satisfy \eqref{inf-sup_discrete}, may be difficult to find. Nevertheless, we can solve \eqref{discrete:variationalSPP} and  avoid building a basis for $\M_h$ by using  an Uzawa type algorithm, e.g., the Uzawa Conjugate Gradient (UCG) algorithm.
\begin{algorithm} (UCG) Algorithm \label{alg:UCG}
\vspace{0.1in}

 {\bf Step 1:} {\bf Choose any}  $p_0 \in \M_h$. {\bf Compute} $u_{1} \in V_h $, $q_1, d_1 \in \M_h$ by 
\[ 
 \begin{aligned}
& a( u_{1}, v_h)& = &\ \<f,v_h\> - b(v_h, p_{0}) &\ &\Forall v_h \in V_h,&\\
&   (q_1, q_h) & = & \ b(u_1 ,q_h)  &\ &\Forall  q_h \in \M_h,& \ \ d_1:=q_1.
\end{aligned}
\]
\vspace{0.1in}

{\bf Step 2:} {\bf For} $j=1,2,\ldots, $ {\bf compute} 
 $h_j, \alpha_j, p_j, u_{j+1}, q_{j+1}, \beta_j, d_{j+1}$ by 
\[ 
 \begin{aligned}
& {\bf (UCG1)} \ \ \ \  & a( h_{j}, v_h) = & - b(v_h, d_j)  && \text{for all} \ v_h \in V_h&\\
& {\bf (UCG\alpha)} \ \ \ \  & \alpha_j =& - \frac{(q_j, q_j)}{b(h_j,q_j)} \\
& {\bf (UCG2)} \ \ \ \  & p_{j} = & \ p_{j-1} + \alpha_j \  d_j \\
& {\bf (UCG3)} \ \ \ \  & u_{j+1} = & \ u_j + \alpha_j\ h_j \\
& {\bf (UCG4)} \ \ \ \  & (q_{j+1}, q_h) = & \ b(u_{j+1} ,q_h) && \text{for all} \ q_h \in \M_h& \\
& {\bf (UCG\beta)} \ \ \ \  & \beta_j=& \ \frac{(q_{j+1}, q_{j+1})}{(q_j,q_j)} \\
& {\bf (UCG6)} \ \ \ \  & d_{j+1}= & \ q_{j+1} +\beta_j d_j. \\
\end{aligned}
\]
\end{algorithm}
Note that the only inversions needed in the algorithm involve the form $a(\cdot,\cdot)$ in \textbf{Step 1} and (\textbf{UCG1}). In operator form, these steps become 
\begin{equation}\label{inversions}
u_1=A_h^{-1}(f_h-B_h^*p_0), \qquad \mathrm{and} \qquad h_j=-A_h^{-1}(B_h^*d_j),
\end{equation}
respectively. In practical  implementations of  Algorithm \ref{alg:UCG}, we  would like  to replace the action of $A_h^{-1}$ with the action of a suitable preconditioner. The properties of the new  preconditioned algorithm  are discussed in the next section. 
The following {\it sharp error estimation} result was proved in \cite{B14}.
\begin{theorem}\label{th:sharp4iterCG} 
If $(u_h, p_h)$ is the discrete solution of \eqref{discrete:variationalSPP} and $(u_{j},p_{j-1})$ is the $j^{th}$ iteration for Algorithm \ref{alg:UCG}, then 
$(u_{j},p_{j-1}) \to (u_{h}, p_h)$ and 

\begin{equation}\label{eq:qInicCG}
\begin{aligned}
\frac{1} {M^2}\, \|q_{j}\|  \leq \|p_{j-1} -p_h \| \leq  \frac{1} {m_h^2}\, \|q_{j}\|, \\
\frac{m_h} {M^2}\, \|q_{j}\|  \leq |u_{j} -u_h | \leq \frac{M} {m_h^2}\,  \|q_{j}\|. 
\end{aligned}
\end{equation} 
%Furthermore, $\|q_j\|\to 0$.
\end{theorem}

\begin{remark}\label{CG-remark}
In particular, Algorithm \ref{alg:UCG} recovers the steps of the conjugate gradient algorithm for solving the Schur complement problem \eqref{Sh}. Hence, the rate of convergence for the iteration error $\|p_{j} -p_h \|_{S_h} $ or  $\|p_{j} -p_h \|$ depends on the condition number of $S_h$, which is $\kappa(S_h)= \frac{M_h^2}{m_h^2}$. 
\end{remark}
%The estimates \eqref{eq:qInicCG} can be used to build adaptive or multilevel algorithms for SPLS discretization. 

%%%%%%%%%%%%%%%%%%%%%%%%%%%%%%%%%%%%%%%5
\section{Preconditioning techniques}
In this section, we develop a general preconditioning framework to approximate the solution of  \eqref{cont_problem} based on \eqref{discrete:variationalSPP}  and elliptic preconditioning of the  operator  associated with the inner product on $V_h$. More precisely, we replace the original form $a(\cdot,\cdot)$ in  \eqref{discrete:variationalSPP}  with a uniformly equivalent form $\tilde{a}(\cdot,\cdot)$ on $V_h$ that leads to an implementably fast operator $\tilde{A}^{-1}_h$. We assume that $V_h\subset V$ and $\M_h\subset Q$ are  finite dimensional approximation spaces satisfying \eqref{inf-sup_discrete} and \eqref{sup-sup_h}.

\subsection{The preconditioned saddle point problem}\label{prec_problem}
First, we introduce a general preconditioner  operator  $P_h:V_h^*\to V_h$ that is  equivalent to  $A_h^{-1}$ in the following sense 
\begin{equation}\label{dual_relation}
\langle g,P_hf\rangle=\langle f,P_hg\rangle \Forall f,g\in V_h^*,
\end{equation}
and
\begin{equation}\label{equiv_forms}
m_1^2|v_h|^2\leq a(P_hA_hv_h,v_h)\leq m_2^2|v_h|^2,
\end{equation}
where the  positive constants $m_1^2,m_2^2$ are the smallest and largest eigenvalues of $P_hA_h$, respectively. 
\begin{remark}
Assumption \eqref{equiv_forms} gives us that the condition number of $P_hA_h$ satisfies
\begin{equation}\label{condition_bound}
\kappa(P_hA_h) = \frac{m_2^2}{m_1^2}.
\end{equation}
\end{remark}

With the  preconditioner $P_h:V_h^*\to V_h$, we  define the form $\tilde{a}:V_h\times V_h\to \mathbb{R}$ by 
\begin{equation}\label{tilde-form-def}
\tilde{a}(u_h,v_h):=a((P_hA_h)^{-1}u_h,v_h) \Forall u_h,v_h\in V_h.
\end{equation}

\begin{prop}\label{symm_PD_prop}
Under assumptions \eqref{dual_relation} and \eqref{equiv_forms}, we have that $\tilde{a}(\cdot, \cdot)$ is symmetric and equivalent with $a(\cdot, \cdot)$ on  $V_h$.
\begin{proof}
For symmetry, it suffices to prove that $P_hA_h$ is symmetric w.r.t. the $a(\cdot,\cdot)$ inner product. From the definition of the operator $A_h$ and \eqref{dual_relation}, we have
\begin{align*}
a(P_hA_hu_h,v_h)=\langle A_h v_h,P_hA_hu_h\rangle&=\langle A_h u_h,P_hA_hv_h\rangle\\
&=a(u_h,P_hA_hv_h).
\end{align*}
For the equivalence, note that \eqref{equiv_forms} and \eqref{tilde-form-def} imply
\begin{equation}\label{norm_equiv_tilde}
\frac{1}{m_2^2}|v_h|^2\leq \tilde{a}(v_h,v_h)\leq \frac{1}{m_1^2}|v_h|^2.
\end{equation}
\end{proof}
\end{prop}
By Proposition \ref{symm_PD_prop}, $\tilde{a}(\cdot,\cdot)$ defines an equivalent inner product on $V_h$. Let $|v_h|_P:=\tilde{a}(v_h,v_h)^{1/2}$ be the norm induced by the inner product  $\tilde{a}(\cdot,\cdot)$  and define the operator $\tilde{A}_h:V_h\to V_h^*$ by
\[
\langle \tilde{A}_h u_h,v_h\rangle:=\tilde{a}(u_h,v_h) \Forall u_h,v_h\in V_h.
\] 
Note that for any $u_h,v_h\in V_h$
\begin{align*}
\langle \tilde{A_h}u_h,v_h\rangle=\tilde{a}(u_h,v_h)&=a((P_hA_h)^{-1}u_h,v_h)\\
&=\langle A_h(P_hA_h)^{-1}u_h,v_h\rangle,
\end{align*}
which implies $\tilde{A}_h=A_h(P_h A_h)^{-1}=P_h^{-1}$. Hence, we can view $\tilde{a}(\cdot,\cdot)$ as a preconditioned version of the form $a(\cdot,\cdot)$. The preconditioned discrete saddle point problem consists of finding $(\tilde{u}_h,\tilde{p}_h)\in V_h\times \M_h$ such that 
\eqref{discrete:variationalSPPprec} holds. To simplify the notation, we will drop the $\tilde{}$ notation from $(\tilde{u}_h,\tilde{p}_h)$. Thus, for the remainder of this paper, the {\it preconditioned saddle point least squares} formulation is: Find $({u}_h,p_h)\in V_h\times \M_h$ such that 
\begin{equation}\label{discrete:variationalSPPpr}
\begin{array}{lclll}
\tilde{a}({u}_h,v_h) & + & b( v_h, p_h) &= \langle f,v_h \rangle &\ \Forall  v_h \in V_h,\\
b({u}_h,q_h) & & & =0   &\  \Forall  q_h \in \M_h.
\end{array}
\end{equation}

Using that $V_h\subset V$ and $\M_h\subset Q$  satisfy \eqref{inf-sup_discrete} and \eqref{sup-sup_h}, we obtain
\begin{equation}\label{inf-sup-tilde}
\tilde{m}_h:=\inf_{p_h\in \M_h}\sup_{v_h\in V_h}\frac{b(v_h,p_h)}{|v_h|_P\,\|p_h\|}\geq m_1 \, m_h >0, 
\end{equation}
and 
\begin{equation}\label{sup-sup-tilde}
\tilde{M}_h:=\sup_{p_h\in \M_h}\sup_{v_h\in V_h}\frac{b(v_h,p_h)}{|v_h|_P\,\|p_h\|}\leq m_2 \, M_h \leq m_2\, M.
\end{equation}
Hence, the {\it preconditioned saddle point least squares} formulation \eqref{discrete:variationalSPPpr} has a unique solution.

The Schur complement associated with problem \eqref{discrete:variationalSPPpr} is 
\[
\tilde{S}_h=B_h\tilde{A}_h^{-1}B_h^*=B_hP_hB_h^*.
\]
 Solving for $p_h$ from   \eqref{discrete:variationalSPPpr}, we obtain 
 \begin{equation} \label{Sh_prec}
\tilde{S}_h \, p_h= B_h (P_hB_h^*) \,p_h = B_h P_h f_h.
\end{equation}
 We call  the component $p_h$  of the solution $(u_h,p_h)$ 
of \eqref{discrete:variationalSPPpr} 
the {\it  (preconditioned) saddle point least squares}  approximation of the solution $p$ of the original mixed prolem \eqref{cont_problem}. To estimate $\|p-p_h\|$ in this case, we will prove the analog to Theorem \ref{th:sharpEE} based on the Xu-Zikatanov argument, see \cite{xu-zikatanov-BBtheory}.
%Using  the extra compatibility condition  \eqref{eq:BBsuf},   a sharp error estimate for $\|p-p_h\|$ was proved in \cite{BQ15}  based on the  Xu-Zikatanov argument, see \cite{xu-zikatanov-BBtheory}. 
\begin{theorem}\label{th:sharpEEpr}  
Let $b:V \times Q \to \R$  satisfy \eqref{inf-sup_a} and \eqref{sup-sup_a} and assume that   ${f}  \in V^*$  is given and satisfies \eqref{comp_condition}. Assume that $V_h \subset V$,  $ {\M}_h \subset  Q$ are  chosen such that the discrete $\inf-\sup$ condition   \eqref{inf-sup_discrete} holds. If  $p$  is the  solution  of \eqref{cont_problem} and  $\left (u_h, p_h \right )$ is the  solution  of \eqref{discrete:variationalSPPpr}, then the following error estimate holds:
\begin{equation}\label{prec_estimate}
\frac{1}{M}\frac{1}{m_2^2}|u_h|\leq \|p-p_h\|\leq \frac{M}{m_h}\frac{m_2}{m_1}\inf_{q_h\in \M_h}\|p-q_h\|.
\end{equation}
\begin{proof}
Define the operator $T_h:Q\to Q$ by $T_hp=p_h$. Note that $T_h$ is linear and idempotent. To show the latter, consider the problem: Find $(u_h^*,p_h^*)\in V_h \times \M_h$ such that 
\begin{equation}\label{discrete:variationalSPPprec2}
\begin{array}{lclll}
\tilde{a}(u_h^*,v_h) & + & b( v_h, p_h^*) &= b(v_h,p_h) &\ \Forall  v_h \in V_h,\\
b(u_h^*,q_h) & & & =0   &\  \Forall  q_h \in \M_h.
\end{array}
\end{equation}
Since $b$ satisfies \eqref{inf-sup_discrete}, we have that \eqref{inf-sup-tilde} is satisfied as described above. Thus, problem \eqref{discrete:variationalSPPprec2} has a unique solution. Since $(u_h^*,p_h^*)=(0,p_h)$ solves the problem, we conclude $T_hp_h=p_h$ which gives us $T_h^2=T_h$. From Kato \cite{kato} and Xu and Zikatanov \cite{xu-zikatanov-BBtheory}, this implies
\[
\|I-T_h\|_{\L(Q,Q)}=\|T_h\|_{\L(Q,Q)}.
\]
Using the above equality, for an arbitrary $q_h\in \M_h$ we have 
\begin{equation}\label{right_ineq}
\|p-q_h\|=\|(I-T_h)p\|=\|(I-T_h)(p-q_h)\|\leq\|T_h\|\,\|p-q_h\|.
\end{equation}
We now estimate $\|T_h\|$. First, define $\tilde{V}_{h,0}^{\perp}$ to be the orthogonal complement of $V_{h,0}$ w.r.t. the $\tilde{a}(\cdot,\cdot)$ inner product. Note that from the first equation of \eqref{discrete:variationalSPPpr} and the fact $p$ solves \eqref{cont_problem} we have that
\begin{equation}\label{first_eq}
b(v_h,p_h)=b(v_h,p)-\tilde{a}(u_h,v_h).
\end{equation}
Also, since $b$ satisfies \eqref{sup-sup_a} we have that \eqref{sup-sup-tilde} holds. Hence, from \eqref{inf-sup-tilde}, \eqref{sup-sup-tilde}, and \eqref{first_eq} we obtain
\begin{align*}
\|T_hp\|\leq \frac{1}{m_h\,m_1}\sup_{v_h\in V_h}\frac{b(v_h,T_hp)}{|v_h|_P}&=\frac{1}{m_h\,m_1}\sup_{v_h\in \tilde{V}_{h,0}^{\perp}}\frac{b(v_h,p_h)}{|v_h|_P}\\
&=\frac{1}{m_h\,m_1}\sup_{v_h\in \tilde{V}_{h,0}^{\perp}}\frac{b(v_h,p)-\tilde{a}(u_h,v_h)}{|v_h|_P}\\
&\leq \frac{Mm_2}{m_hm_1}\|p\|. \stepcounter{equation}\tag{\theequation}\label{right_ineq2}
\end{align*}
The right inequality now follows from \eqref{right_ineq} and \eqref{right_ineq2}.
For the left inequality, note that 
\begin{align*}
|u_h|_P=\sup_{v_h\in V_h}\frac{\tilde{a}(u_h,v_h)}{|v_h|_P}&=\sup_{v_h\in V_h}\frac{b(v_h,p-p_h)}{|v_h|_P}\leq M\,m_2\|p-p_h\|,
\end{align*}
and 
\[
|u_h|\leq m_2|u_h|_P.
\]
\end{proof}
\end{theorem}

\subsection{An iterative solver for the preconditioned variational formulation}\label{iter_prec_problem}
We use a modified version of Algorithm \ref{alg:UCG} to solve \eqref{discrete:variationalSPPpr} by replacing the form $a(\cdot,\cdot)$ by $\tilde{a}(\cdot,\cdot)$ in \textbf{Step 1} and (\textbf{UCG1}). With this modification, we obtain the following (Uzawa) Preconditioned Conjugate Gradient (PCG) algorithm for mixed methods.
\begin{algorithm} (PCG) Algorithm for Mixed Methods\label{alg:PUCG}
\vspace{0.1in}

 {\bf Step 1:} {\bf Choose any} $p_0 \in \M_h$. {\bf Compute} $u_{1} \in V_h $, $q_1, d_1 \in \M_h$ by 
\[ 
 \begin{aligned}
& u_1 &=&P_h(f_h-B_h^*p_0) \\
& q_1& =& B_h u_1, \ \ d_1:=q_1.
\end{aligned}
\]
\vspace{0.1in}

{\bf Step 2:} {\bf For} $j=1,2,\ldots, $ {\bf compute} 
 $h_j, \alpha_j, p_j, u_{j+1}, q_{j+1}, \beta_j, d_{j+1}$ by 
\[ 
 \begin{aligned}
& {\bf (PCG1)} \ \ \ \  & h_j =& - P_h(B_h^*d_j)  && &\\
& {\bf (PCG\alpha)} \ \ \ \  & \alpha_j =& - \frac{(q_j, q_j)}{b(h_j,q_j)} \\
& {\bf (PCG2)} \ \ \ \  & p_{j} = & \ p_{j-1} + \alpha_j \  d_j \\
& {\bf (PCG3)} \ \ \ \  & u_{j+1} = & \ u_j + \alpha_j\ h_j \\
& {\bf (PCG4)} \ \ \ \  & q_{j+1} = & B_h u_{j+1} ,& & & \\
& {\bf (PCG\beta)} \ \ \ \  & \beta_j=& \ \frac{(q_{j+1}, q_{j+1})}{(q_j,q_j)} \\
& {\bf (PCG6)} \ \ \ \  & d_{j+1}= & \ q_{j+1} +\beta_j d_j. \\
\end{aligned}
\]
\end{algorithm}
Note that only the actions of $P_h$, $B_h$,  and $B^*_h$ are needed in the above algorithm. For any preconditioner $P_h$ and trial space $\M_h$ that is not defined via a global projection, these actions  do not involve inversion processes, see Section \ref{multilevel-example} for the case $P_h$ -an additive multilevel Schwarz preconditioner. Similar to the remark in Section \ref{iter_discrete_problem}, we have the following:
\begin{remark}\label{UPCG-remark}
 Algorithm \ref{alg:PUCG} recovers in particular the steps of the conjugate gradient algorithm for solving the problem \eqref{Sh_prec}. Hence, the rate of convergence for $\|p_{j} -p_h \|_{\tilde{S}_h} $ or  $\|p_{j} -p_h \|$ depends on the condition number of $\tilde{S}_h$, which is $\kappa(\tilde{S}_h)= \frac{\tilde{M}_h^2}{\tilde{m}_h^2}$. 
\end{remark}
The following result is analogous to Theorem \ref{th:sharp4iterCG}.
\begin{theorem}\label{th:sharp4iterPCG} 
If $(u_h, p_h)$ is the discrete solution of \eqref{discrete:variationalSPPpr} and $(u_{j},p_{j-1})$ is the $j^{th}$ iteration for Algorithm \ref{alg:PUCG}, then 
$(u_{j},p_{j-1}) \to (u_{h}, p_h)$ and 

\begin{equation}\label{eq:qInicPCG}
\begin{aligned}
\frac{1}{M^2}\frac{1}{m_2^2}\, \|q_{j}\|  \leq \|p_{j-1} -p_h \| \leq  \frac{1}{m_h^2}\frac{1}{m_1^2}\, \|q_{j}\|, \\
\frac{m_h} {M^2}\frac{m_1^2}{m_2^2}\, \|q_{j}\|  \leq |u_{j} -u_h | \leq \frac{M} {m_h^2}\frac{m_2^2}{m_1^2}\,  \|q_{j}\|. 
\end{aligned}
\end{equation} 
%Furthermore, $\|q_j\|\to 0$.
\begin{proof}
By induction over $j$, we have that
\[
\tilde{a}(u_{j},v_h)+b(v_h,p_{j-1})=\langle f,v_h\rangle \Forall v_h\in V_h.
\]
Combining this with the first equation of \eqref{discrete:variationalSPPpr} gives us
\begin{equation}\label{tilde_eq1}
\tilde{a}(u_{j}-u_h,v_h)=b(v_h,p_h-p_{j-1}) \Forall v_h\in V_h.
\end{equation} 
Note that $\sigma(\tilde{S}_h)\subset [\tilde{m}_h^2,\tilde{M}_h^2]$. Hence,
\begin{equation}\label{schur_est_tilde}
\tilde{m}_h\|q_h\|=(\tilde{S}_hq_h,q_h)^{1/2}\leq \tilde{M}_h\|q_h\| \Forall q_h\in \M_h.
\end{equation}
By substituting $v_h=\tilde{A}_h^{-1}B_h^*(p_h-p_{j-1})$ into \eqref{tilde_eq1},
\[
|u_j-u_h|_P^2=(\tilde{S}_h(p_h-p_{j-1}),p_h-p_{j-1})=\|p_h-p_{j-1}\|_{\tilde{S}_h}^2.
\]
The above equality, \eqref{norm_equiv_tilde}, and \eqref{schur_est_tilde} gives us that
\begin{equation}\label{est1}
m_1\tilde{m}_h\|p_h-p_{j-1}\|\leq |u_j-u_h|\leq m_2\tilde{M}_h\|p_h-p_{j-1}\|.
\end{equation}
From (\textbf{PCG4}), the second equation of \eqref{discrete:variationalSPPpr}, and \eqref{tilde_eq1} we have that
\[
q_{j}=B_hu_{j}=B_h(u_{j}-u_h)=\tilde{S}_h(p_h-p_{j-1}).
\]
Thus,
\begin{equation}\label{est2}
\tilde{m}_h^2\|p_h-p_{j-1}\|\leq \|\tilde{S}_h(p_h-p_{j-1})\|=\|q_j\|\leq \tilde{M}_h^2\|p_h-p_{j-1}\|.
\end{equation}
The inequalities \eqref{eq:qInicPCG} follow from \eqref{est1}, \eqref{est2}, and the fact that $\tilde{m}_h\geq m_hm_1$ and $\tilde{M}_h\leq Mm_2$. From  Remark \ref{UPCG-remark} and the standard estimate for the convergence rate of the conjugate gradient algorithm, \cite{braess,hestenes}, we have that
\begin{equation}\label{PCGrate}
\|p_h-p_j\|_{\tilde{S}_h}\leq 2\left(\frac{\tilde{M}_h-\tilde{m}_h}{\tilde{M}_h+\tilde{m}_h}\right)^j\|p_h-p_0\|_{\tilde{S}_h}. 
\end{equation}
Hence, $p_{j}\to p_h$. From \eqref{eq:qInicPCG}, we conclude that $u_{j} \to u_h$ as well. 
\end{proof}
\end{theorem}
The following estimates are a direct consequence of \eqref{inf-sup-tilde}, \eqref{sup-sup-tilde}, \eqref{PCGrate}, and the formula $\kappa(\tilde{S}_h)=\frac{\tilde{M}_h^2}{\tilde{m}_h^2}$.
\begin{prop} The condition number of the  Schur complement \\ $\tilde{S}_h=B_hP_hB_h^*$ satisfies
\begin{equation}\label{condition-number-estimate}
\kappa(\tilde{S}_h)\leq  \frac{M_h^2}{m_h^2}\frac{m_2^2}{m_1^2}=\kappa(S_h)\cdot \kappa(P_hA_h). 
\end{equation}
Consequently,  the convergence rate   $\rho_h$  for $\|p_{j} -p_h \|_{\tilde{S}_h} $ satisfies
\[
\rho_h\leq \frac{\frac{M_h}{m_h} \frac{m_2}{m_1} -1} {\frac{M_h}{m_h} \frac{m_2}{m_1} +1}.
\]
\end{prop}

\begin{remark}
We can relate our preconditioned SPLS discretization method  for  solving the general mixed problem  \eqref {cont_problem} with the {\it Bramble-Pasciak least squares} approach presented in \cite{BPls03}. In our notation, the  Bramble-Pasciak least squares discretization can be formulated as: Find $p_h\in \M_h$ such that
\[
b(A_h^{-1}B_h^*q_h,p_h)=\langle f_h, A_h^{-1}B_h^*q_h\rangle=b(A_h^{-1}f_h,q_h) \Forall q_h\in \M_h. 
\]
With a suitable preconditioner $P_h$ replacing $A_h^{-1}$, the problem becomes: Find $p_h\in \M_h$ such that
\begin{equation}\label{LS_prec_problem}
b(P_hB_h^*q_h,p_h)=b(P_hf_h,q_h) \Forall q_h\in \M_h.
\end{equation}
We shall note that \eqref{LS_prec_problem}   is equivalent to our  Schur complement problem \eqref{Sh_prec}. While we arrive at essentially the same normal equation for  solving \eqref{discrete_var_form}, our saddle point approach is more direct and allows 
sharp error estimates for the error $\|p-p_h\|$.  The two approaches are also essentially different in the way the trial spaces are chosen, see Section \ref{discrete_spaces} for our choices of trial spaces.  In  \cite{BPls03}, to iteratively solve \eqref{LS_prec_problem},  bases for both the test and trial spaces are needed. In contrast, we solve the coupled preconditioned saddle point problem \eqref{discrete:variationalSPPpr} using Algorithm \ref{alg:PUCG} which avoids the need of a basis for the trial space. %The $p_h$ component of the solution to \eqref{discrete:variationalSPPpr} is in fact the solution to \eqref{LS_prec_problem}.
\end{remark}

\subsection{An example of  a preconditioner.}\label{multilevel-example}
In order to illustrate the  applicability  of the theory presented thus far, we consider the case when $P_h$ is given by the additive multilevel Schwarz or BPX preconditioner, see \cite{BPX90, BrambleZhang,Zhang92}. Assume that we have a nested sequence of approximation spaces $V_1\subset V_2\subset\cdots \subset V_J=V_h$ and let $\{\phi_1^k,\phi_2^k,\dots,\phi_{n_k}^k\}$ be a basis for $V_k$. For $f_h\in V_h^*$, the action of $P_h$ is given by
\[
P_hf_h=\sum_{k=1}^J\sum_{i=1}^{n_k}\frac{\langle f_h,\phi_i^k\rangle}{a(\phi_i^k,\phi_i^k)}\phi_i^k.
\]

It is known that for $V=H_0^1(\Omega)$ and a nested sequence $\{V_k\}$ of piecewise linear functions that,  under standard  mesh uniformity conditions,  $P_h$ is a preconditioner for $A_h$ satisfying \eqref{dual_relation} and \eqref{equiv_forms}, see \cite{BPX90,JungBPX02,XuSIAMReview,Xu-Qin94,Zhang92}.

In this case, the first equation in \textbf{(Step 1)} of Algorithm \ref{alg:PUCG} becomes
\[
u_1=P_h(f_h-B_h^*p_0)=\sum_{k=1}^J\sum_{i=1}^{n_k}\frac{\langle f_h,\phi_i^k\rangle-b(\phi_i^k,p_0)}{a(\phi_i^k,\phi_i^k)}\phi_i^k.
\]
Furthermore, the iterates for $h_j$ in \textbf{(PCG1)} are given by
\[
h_j=-\sum_{k=1}^J\sum_{i=1}^{n_k}\frac{b(\phi_i^k,d_j)}{a(\phi_i^k,\phi_i^k)}\phi_i^k,
\]
which implies that
\[
b(h_j,q_j)=-\sum_{k=1}^J\sum_{i=1}^{n_k}\frac{b(\phi_i^k,d_j)b(\phi_i^k,q_j)}{a(\phi_i^k,\phi_i^k)},
\]
in \textbf{(PCG$\alpha$)}. Thus, the implementation of Algorithm \ref{alg:PUCG} does not involve matrix inversion. Certainly, any elliptic preconditioner, including the standard multigrid ones, can be used for $P_h$.  We decided to show  details of a general   additive multilevel Schwarz  (or BPX) preconditioner to emphasize the  simplicity of implementation when dealing with mixed methods preconditioning. More details on implementing the matrix action of multilevel  preconditioners (including BPX)  can be found in \cite{Xu-Qin94}.

%%%%%%%%%%%%%%%%%%%%%%%%%%%%%
\section{Discrete spaces that satisfy an $\inf-\sup$ condition}\label{discrete_spaces}
In this section, we describe two pairs of discrete spaces, introduced in \cite{BQ15}, which satisfy the discrete $\inf-\sup$ condition \eqref{inf-sup_discrete} in the general abstract framework of Section \ref{abstract}. In light of \eqref{condition-number-estimate}, we would like to provide families of spaces $\{(V_h,\M_h)\}$ such that $\kappa(S_h)$ is small. Let $V_h\subset V$ be a finite element test space and assume the action of $\C^{-1}$, where $\C$  was defined in Section \ref{abstract},  is easy to obtain at the continuous level. 

\subsection{No projection trial space.}\label{no_proj} 
The first choice defines $\M_h\subset Q$ by
\[
\M_h:=\C^{-1}BV_h.
\]
In this case, $V_{h,0}\subset V_0$ and a discrete $\inf-\sup$ condition holds. Indeed, for a generic $p_h =\C^{-1}Bw_h\in \M_h$ where $w_h\in V_{h,0}^{\perp}$, we have 
\begin{align*}
m_{h,0}&:=\inf_{p_h \in \M_h}\sup_{v_h\in V_h} \frac{b(v_h,p_h)}{|v_h|\  \|p_h\|}=\inf_{w_h \in V^{\perp}_{h,0}}\sup_{v_h\in V_h} \frac{(\C^{-1}Bv_h,\C^{-1}Bw_h)}{|v_h|\,\|\C^{-1}Bw_h\|}\\ \stepcounter{equation}\tag{\theequation}\label{inf-supNoProj}
&\geq \inf_{w_h \in V^{\perp}_{h,0}}\frac{\|\C^{-1}Bw_h\|^2}{|w_h|\,\|\C^{-1}Bw_h\|} 
=\inf_{w_h \in V^{\perp}_{h,0}}\frac{\|\C^{-1}Bw_h\|}{|w_h|}
>0.
\end{align*}
Hence, by Remark \ref{well-posed_prop} we have that \eqref{discrete_var_form} has a unique solution $p_h\in \M_h$ and $(u_h=0,p_h)$ solves \eqref{discrete:variationalSPP}. In this case, $p_h$ is an optimal approximation to the solution $p$ of \eqref{cont_problem}. Indeed, for any $v_h\in V_h$
\begin{align*}
0=b(v_h,p-p_h)&=\langle Bv_h,p-p_h\rangle\\
&=(\C^{-1}Bv_h,p-p_h). 
\end{align*}
Thus, $p_h$ is the orthogonal projection of $p$ onto $\M_h$ which implies
\begin{equation}\label{optimal_est}
\|p-p_h\|=\inf_{q_h\in \M_h}\|p-q_h\|.
\end{equation}
While \eqref{optimal_est} gives optimal approximation error, to efficiently approximate $p_h$ using Algorithm \ref{alg:UCG} or \ref{alg:PUCG}, it requires spaces $\{(V_h,\M_h)\}$ for which $\kappa(S_h)=\frac{M_h^2}{m_h^2}$ is small or independent of $h$.
\subsection{Projection type trial space.}\label{proj}
The second choice defines $\M_h\subset Q$ by
\[
\M_h:=R_h\C^{-1}BV_h,
\]
where $R_h:Q\to \tilde{M}_h$ is defined by
\begin{equation}\label{R_h_def}
(R_hp,q_h)_h:=(p,q_h) \Forall q_h\in \tilde{M}_h,
\end{equation}
and $\tilde{M}_h$ is a finite dimensional subspace of $Q$ equipped with the inner product $(\cdot, \cdot)_h$. 
\begin{remark}
If the $(\cdot, \cdot)_h$ inner product coincides with the inner product on $Q$, then by definition $R_h$ is the orthogonal projection onto $\tilde{M}_h$. 
\end{remark}
In general, the inner product on $\tilde{M}_h$ could be different from the inner product on $Q$, but we assume that $(\cdot,\cdot)_h$ induces an equivalent norm independent of $h$. The following proposition provides a sufficient condition on $R_h$ which implies well-posedness of problems \eqref{discrete:variationalSPP} and \eqref{discrete_var_form} and relates the stability of the family of spaces $\{(V_h,R_h\C^{-1}BV_h)\}$ with the stability of the family of spaces $\{(V_h,\C^{-1}BV_h)\}$ defined in Section \ref{no_proj}.
\begin{prop}\label{R_h_trial_prop}
Assume that
\begin{equation}\label{R_h_coercive}
\|R_hq_h\|_h\geq \tilde{c}\|q_h\| \Forall q_h\in \C^{-1}BV_h,
\end{equation}
with a constant $\tilde{c}$ independent of $h$. Then $V_{h,0}\subset V_0$. Furthermore, the stability of the family $\{(V_h,\C^{-1}BV_h)\}$ implies the stability of the family $\{V_h,R_h\C^{-1}BV_h)\}$.
\begin{proof}
Let $v_h \in V_{h,0}$. Then, for any $p_h \in \M_h$,  
\[
0=b(v_h,p_h)=(\C^{-1}Bv_h,p_h)=(R_h\C^{-1}Bv_h,p_h)_h.
\]
Taking $p_h=R_h\C^{-1}Bv_h$ gives us $\|R_h\C^{-1}Bv_h\|_h=0$ and the inclusion $V_{h,0} \subset V_0$ follows from \eqref{R_h_coercive}. For the stability, note that for a generic function $p_h =R_h\C^{-1}Bw_h\in \M_h$, where $w_h \in V^{\perp}_{h,0}$, we have
\begin{align*}
m_h=\inf_{p_h \in \M_h}\sup_{v_h\in V_h} \frac{b(v_h,p_h)}{|v_h|\  \|p_h\|_h}&= \inf_{w_h \in V^{\perp}_{h,0}}\sup_{v_h\in V_h} \frac{(\C^{-1}Bv_h,R_h\C^{-1}Bw_h)}{|v_h|\  \|R_h\C^{-1}Bw_h\|_h}\\
&=\inf_{w_h \in V^{\perp}_{h,0}}\sup_{v_h\in V_h} \frac{(R_h\C^{-1}Bv_h,R_h\C^{-1}Bw_h)_h}{|v_h|\  \|R_h\C^{-1}Bw_h\|_h}\\
&\geq \inf_{w_h \in V^{\perp}_{h,0}} \frac{\|R_h\C^{-1}Bw_h\|_h^2}{|w_h|\  \|R_h\C^{-1}Bw_h\|_h}\\
&\geq \tilde{c} \inf_{w_h \in V^{\perp}_{h,0}} \frac{\|\C^{-1}Bw_h\|}{|w_h|}
=\tilde{c}\, m_{h,0},
\end{align*}
where $m_{h,0}$ is defined in \eqref{inf-supNoProj}.
\end{proof}
\end{prop}
In this case, we have that $p_h$ is a quasi-optimal approximation of the solution $p$ of \eqref{cont_problem} by Theorem \ref{th:sharpEEpr}. 

The benefit of using the projection type trial space is that it could lead to a better approximation of the continuous solution $p$. Indeed, for the case when preconditioning is not used, super-convergence of $\|p-p_h\|$ is observed, see \cite{BJQS17, BQ15,BQ17}. Using uniform preconditioners and Theorems \ref{th:sharpEEpr} and \ref{th:sharp4iterPCG}, we expect 
the same order of  super-convergence for  $\|p-p_h\|$.  
\begin{remark}
With this choice of trial space, Algorithms \ref{alg:UCG} and \ref{alg:PUCG} need to be modified to account for the $(\cdot,\cdot)_h$ inner product on $\M_h\subset\tilde{M}_h$. This modification is nothing more than replacing the $(\cdot,\cdot)$ inner product with the $(\cdot,\cdot)_h$ inner product where it appears in the algorithms.
\end{remark}

%%%%%%%%%%%%%%%%%%%%%%%%%%%%%%%%%
\section{Conclusion}
We presented a general preconditioning approach to mixed variational formulations of the form \eqref{cont_problem} that relies on 
the classical theory of symmetric saddle point problems and on  the theory of preconditioning symmetric positive definite operators. 
First, a discrete saddle point variational formulation \eqref{discrete:variationalSPP}, that approximates the solution of the original mixed problem in a  least squares sense, is considered. In this formulation, the inner product $a(\cdot,\cdot)$ is replaced by an equivalent bilinear form that give rise to efficient  elliptic inversion or preconditioning.  An Uzawa preconditioned conjugate gradient algorithm for solving the new saddle point system was proposed that requires bases only for the space $V_h$ and avoids costly inversion processes. Due to the saddle point interpretation of the preconditioned system, we were able to prove sharp approximability properties for the discretization and iteration errors and were able to provide practical estimates for the rate of convergence of the final preconditioned conjugate algorithm. Using a common test space, two choices of compatible discrete spaces were given.

We plan to apply preconditioned saddle point least squares discretization to first order systems of PDEs as well as second order elliptic problems. In addition, we plan to combine this approach with multilevel and adaptive techniques.
\bibliography{bacutaBib}

\def\cprime{$'$} \def\ocirc#1{\ifmmode\setbox0=\hbox{$#1$}\dimen0=\ht0
  \advance\dimen0 by1pt\rlap{\hbox to\wd0{\hss\raise\dimen0
  \hbox{\hskip.2em$\scriptscriptstyle\circ$}\hss}}#1\else {\accent"17 #1}\fi}
  \def\cprime{$'$} \def\ocirc#1{\ifmmode\setbox0=\hbox{$#1$}\dimen0=\ht0
  \advance\dimen0 by1pt\rlap{\hbox to\wd0{\hss\raise\dimen0
  \hbox{\hskip.2em$\scriptscriptstyle\circ$}\hss}}#1\else {\accent"17 #1}\fi}
\begin{thebibliography}{10}

\bibitem{A-B}
A.~Aziz and I.~Babu\v{s}ka.
\newblock Survey lectures on mathematical foundations of the finite element
  method.
\newblock {\em The Mathematical Foundations of the Finite Element Method with
  Applications to Partial Differential Equations, A. Aziz, editor}, 1972.

\bibitem{B09}
C.~Bacuta.
\newblock Schur complements on {H}ilbert spaces and saddle point systems.
\newblock {\em J. Comput. Appl. Math.}, 225(2):581--593, 2009.

\bibitem{B14}
C.~Bacuta.
\newblock Cascadic multilevel algorithms for symmetric saddle point systems.
\newblock {\em Comput. Math. Appl.}, 67(10):1905--1913, 2014.

\bibitem{BJQS17}
C.~Bacuta, J.~Jacavage, K.~Qirko, and F.J. Sayas.
\newblock Saddle point least squares iterative solvers for the time harmonic
  {M}axwell equations.
\newblock {\em Comput. Math. Appl.}, 74(11):2915--2928, 2017.

\bibitem{BM12}
C.~Bacuta and P.~Monk.
\newblock Multilevel discretization of symmetric saddle point systems without
  the discrete {LBB} condition.
\newblock {\em Appl. Numer. Math.}, 62(6):667--681, 2012.

\bibitem{BQ15}
C.~Bacuta and K.~Qirko.
\newblock A saddle point least squares approach to mixed methods.
\newblock {\em Comput. Math. Appl.}, 70(12):2920--2932, 2015.

\bibitem{BQ17}
C.~Bacuta and K.~Qirko.
\newblock A saddle point least squares approach for primal mixed formulations
  of second order {PDE}s.
\newblock {\em Comput. Math. Appl.}, 73(2):173--186, 2017.

\bibitem{braess}
D.~Braess.
\newblock {\em Finite Elements. Theory, Fast Solvers, and Applications in Solid
  Mechanics}.
\newblock Cambridge University Press, Cambridge, 1997.

\bibitem{BPls03}
J.H. Bramble and J.E. Pasciak.
\newblock A new approximation technique for div-curl systems.
\newblock {\em Math. Comp.}, 73:1739--1762, 2004.

\bibitem{BPX90}
J.H. Bramble, J.E. Pasciak, and J.~Xu.
\newblock Parallel multilevel preconditioners.
\newblock {\em Math. Comp.}, 55(191):1--22, 1990.

\bibitem{BrambleZhang}
J.H. Bramble and X.~Zhang.
\newblock The analysis of multigrid methods.
\newblock In {\em Handbook of numerical analysis, Vol. VII}, pages 173--415.
  North-Holland, Amsterdam, 2000.

\bibitem{Dahmen-Welper-Cohen12}
A.~Cohen, W.~Dahmen, and G.~Welper.
\newblock Adaptivity and variational stabilization for convection-diffusion
  equations.
\newblock {\em ESAIM Math. Model. Numer. Anal.}, 46(5):1247--1273, 2012.

\bibitem{demkowicz-gopalakrishnanDPG10}
L.~Demkowicz and J.~Gopalakrishnan.
\newblock A class of discontinuous {P}etrov-{G}alerkin methods. {P}art {I}: the
  transport equation.
\newblock {\em Comput. Methods Appl. Mech. Engrg.}, 199(23-24):1558--1572,
  2010.

\bibitem{demkowicz-gopalakrishnanDPG13}
L.~Demkowicz and J.~Gopalakrishnan.
\newblock A primal {DPG} method without a first-order reformulation.
\newblock {\em Comput. Math. Appl.}, 66(6):1058--1064, 2013.

\bibitem{hestenes}
M.R. Hestenes and E.~Stiefel.
\newblock Methods of conjugate gradients for solving linear systems.
\newblock {\em J. Research Nat. Bur. Standards}, 49:409--436 (1953), 1952.

\bibitem{JungBPX02}
M.~Jung, S.~Nicaise, and J.~Tabka.
\newblock Some multilevel methods on graded meshes.
\newblock {\em Journal of Computational and Applied Mathematics}, 138(1):151 --
  171, 2002.

\bibitem{kato}
T.~Kato.
\newblock Estimation of iterated matrices, with application to the {V}on
  {N}eumann condition.
\newblock {\em Numer. Math.}, 2:22--29, 1960.

\bibitem{XuSIAMReview}
J.~Xu.
\newblock Iterative methods by space decomposition and subspace correction.
\newblock {\em SIAM Review}, 34:581--613, 1992.

\bibitem{Xu-Qin94}
J.~Xu and J.~Qin.
\newblock Some remarks on a multigrid preconditioner.
\newblock {\em SIAM Journal on Scientific Computing}, 15(1):172--184, 1994.

\bibitem{xu-zikatanov-BBtheory}
J.~Xu and L.~Zikatanov.
\newblock Some observations on {B}abu\v ska and {B}rezzi theories.
\newblock {\em Numer. Math.}, 94(1):195--202, 2003.

\bibitem{Zhang92}
X.~Zhang.
\newblock Multilevel schwarz methods.
\newblock {\em Numerische Mathematik}, 63(1):521--539, Dec 1992.

\end{thebibliography}
\bibliographystyle{plain}

\end{document}